# Approximate Functional Equation for the Product of Functions and Divisor Problem

V.V. Rane
Department of Mathematics
The Institute of Science
15 , Madame Cama Road
Mumbai-400 032
India .
v_v_rane@yahoo.co.in

**Abstract** : For functions defined via Dirichlet series/generalized Dirichlet series in some half planes of the complex plane , we give a new simple elementary approach to obtain an Approximate Functional Equation (AFE for short) for the product of functions with explicit remainder term , provided the given functions have AFEs in the same region . We do this , using a new highly generalized identity , following Motohashi's approach based on a generalization of Dirichlet's device . We use our method to give an AFE for the product of two Dirichlet L-series corresponding to primitive Dirichlet characters , having explicit expression for the remainder term , with hitherto unknown estimate for it . We also consider AFEs of functions on real line and state AFEs of square-roots of Riemann zeta function and Dirichlet L-series in the critical strip .

# Approximate Functional Equation for the Product of Functions and Divisor Problem

V.V. Rane
Department of Mathematics
The Institute of Science
15 , Madame Cama Road
Mumbai-400 032
India .
v_v_rane@yahoo.co.in

**Introduction** : For a complex variable $s = \sigma + it$ , where $\sigma, t$ are real , let

$$f_1(s) = \sum_{n \geq 1} a_1(n) n^{-s}; \quad g_1(s) = \sum_{n \geq 1} b_1(n) n^{-s};$$

$$f_2(s) = \sum_{n \geq 1} a_2(n) n^{-s}; \quad g_2(s) = \sum_{n \geq 1} b_2(n) n^{-s};$$

be given Dirichlet Series , convergent in some half planes of the complex plane . Thus $f_1(s) f_2(s) = \sum_{n \geq 1} A(n) n^{-s}$ and $g_1(s) g_2(s) = \sum_{n \geq 1} B(n) n^{-s}$ in some half planes of the complex plane , where $A(n) = \sum_{uv=n} a_1(u) a_2(v)$ and

$$B(n) = \sum_{uv=n} b_1(u) b_2(v) .$$

In a certain region of the complex plane , let there be relations

$$f_1(s) = \sum_{n \leq x_1} a_1(n) n^{-s} + \psi_1(s) \sum_{n \leq y_1} b_1(n) n^{s-\delta} + E_1(s, x_1); \text{ and}$$

$$f_2(s) = \sum_{n \leq x_2} a_2(n) n^{-s} + \psi_2(s) \sum_{n \leq y_2} b_2(n) n^{s-\delta} + E_2(s, x_2);$$

where $\psi_1(s), \psi_2(s)$ are given functions and $\delta$ is a fixed complex number and

$x_1 y_1 = T_1$ and $x_2 y_2 = T_2$; $x_1, x_2 > 0$. That is, $x_1 > 0, y_1 = \frac{T_1}{x_1}$; and $x_2 > 0, y_2 = \frac{T_2}{x_2}$.

Thus $y_1 \geq 0, y_2 \geq 0$ with $T_1, T_2$ depending upon $s = \sigma + it$ only, and possibly some other parameters. Here $T_1$ is independent of $x_1$; and $T_2$ independent of $x_2$. The motivation for considering these relations comes from the Approximate Functional Equation (AFE) for Riemann Zeta Function $\zeta(s)$ in the critical strip namely,

$$\zeta(s) = \sum_{n \leq x} n^{-s} + \chi(s) \sum_{n \leq y} n^{s-1} + E(s,x)$$ with $x > 0, y = \frac{|t|}{2\pi x}$, where $\chi(s)$ is

defined by the functional equation $\zeta(s) = \chi(s) \zeta(1-s)$. We also know that

$E(s,x) \ll x^{-\sigma} + \frac{x^{1-\sigma}}{\sqrt{t}}$, provided $x \gg 1, y \gg 1$ with $\ll$ constants absolute.

Following the analogy with the case of Riemann zeta function $\zeta(s)$, we shall call the above two relations Approximate Functional Equations (AFEs) for $f_1(s)$ and $f_2(s)$ respectively, though formally neither $f_1(s)$ nor $f_2(s)$ may have any functional equation at all. The object of this paper is to give an AFE for the product $f_1(s) \cdot f_2(s)$ in the same region by elementary method with an explicit expression for its remainder term $E_{1,2}(s, x_1 x_2)$ and to give an estimate for $E_{1,2}(s, x_1 x_2)$, provided we have estimates for both $E_1(s, x_1)$ and $E_2(s, x_2)$. Besides being the elementary method and hence transparent, our approach is significant in the sense that it does away with the impression that

a) the derivation of AFE-type relation for the product of functions is highly

dependent upon the functional equation satisfied by each function.

b) the derivation of AFE for the product depends upon the peculiar nature of the functional equation satisfied by each function.

c) the AFE for the product can be obtained, only by analytic methods comprising complex contour integration.

Our method enables one to write down the AFE, not only of the r-th power of $\zeta(s)$ or $L(s,\chi)$ but also, of the r-th root of $\zeta(s)$ or $L(s,\chi)$, where $r \geq 1$ is an integer. Also our approach brings out new facts. (See our corollaries below).

As regards the AFE for $\zeta(s)$ as stated above, one particular case of it namely, the AFE of $\zeta(s)$ for real $s$ with $0 < s \neq 1$ need be stated. Note that if $0 < s$, then $t = 0$ and consequently, $y = 0$. Thus the AFE of $\zeta(s)$ reads as $\zeta(s) = \sum_{n \leq x} n^{-s} + E(s,x)$ for $x > 0$ arbitrary and $0 < s \neq 1$. In fact, it can be seen that $E(s,x) = \frac{x^{1-s}}{s-1} + x^{-s}\psi(x) - s\int_{x}^{\infty} \frac{\psi(u)}{u^{s+1}} du$, where $\psi(u) = u - [u] - \frac{1}{2}$. Here $[u]$ denotes integral part of real variable $u$. Thus we have

$E(s,x) = \frac{x^{1-s}}{s-1} + x^{-s}\psi(x) + 0(x^{-s}) = \frac{x^{1-s}}{s-1} + 0(x^{-s})$ for $0 < s \neq 1$, where the

0-constant is absolute. Thus, one can introduce the concept of the AFE of $\zeta(s)$ for positive real $s \neq 1$. Following our Theorem, the AFE of $\zeta(s)$ for positive $s$ will result in the AFE of $\zeta^r(s)$ for $s$ positive, where $r \geq 1$ is an integer.

More generally for $0 < \alpha \leq 1$, if $\zeta(s,\alpha)$ is the Hurwitz zeta function

: 4 :

defined by $\zeta(s,\alpha) = \sum_{n \geq 0}(n+\alpha)^{-s}$ for $\sigma > 1$; and its analytic continuation,

then we have an AFE for $\zeta(s,\alpha)$ given by

$$\zeta(s,\alpha) = \sum_{0 \leq n \leq x-\alpha}(n+\alpha)^{-s} + 2^s \pi^{s-1}\Gamma(1-s)\sum_{n \leq y}\sin\left(\frac{\pi s}{2} + 2\pi n\alpha\right)n^{s-1} + E(s,x,\alpha),$$

where $x > 0$ is arbitrary and $y = \frac{|t|}{2\pi x}$. For $x \gg 1$ and $y \gg 1$, we have

$$E(s,x,\alpha) \ll x^{-\sigma} \cdot \log(y+2) + \frac{x^{1-\sigma}}{\sqrt{t}}$$ uniformly in $0 < \alpha \leq 1$ and for $0 < \sigma < 1$.

(Note that $\zeta(s,\alpha) = 2^s \pi^{s-1} \Gamma(1-s)\sum_{n \geq 1}\sin\left(\frac{\pi s}{2} + 2\pi n\alpha\right)n^{s-1}$ for $\sigma < 0$ ).

Here $\ll$ constant is absolute and $\Gamma(s)$ denotes the gamma function. This can be easily checked from author [3]. Also see author [4]. Consequently in particular, for real $s$ this results in the AFE of $\zeta(s,\alpha)$ in the form

$$\zeta(s,\alpha) = \sum_{0 \leq n \leq x-\alpha}(n+\alpha)^{-s} + E(s,x,\alpha),$$ where $x > 0$ is arbitrary and $0 < s < \infty$

and $s \neq 1$. We can see that

$$E(s,x,\alpha) = \frac{x^{1-s}}{s-1} + \psi(x-\alpha)x^{-s} - s\int_x^\infty \frac{\psi(u-\alpha)}{u^{s+1}}du, \text{ where } \psi(u) = u - [u] - \frac{1}{2}.$$

Thus for $s > 0$, $E(s,x,\alpha) = \frac{x^{1-s}}{s-1} + O(x^{-s})$ as in the case of $\zeta(s)$. For a given integer $q \geq 1$ and for an integer a with $1 \leq a \leq q$, let $z(s,a,q) = q^{-s}\zeta(s,\frac{a}{q})$.

Note that if $\chi(\bmod q)$ is a Dirichlet character and if $L(s,\chi) = \sum_{n \geq 1}\chi(n)n^{-s}$ for

$\sigma > 1$; and its analytic continuation, then $L(s,\chi) = \sum_{a=1}^{q}\chi(a)z(s,a,q)$. From

the AFE of $\zeta(s,\alpha)$, we get the AFE of $z(s,a,q)$ in the form

$$z(s,a,q) = \sum_{n \leq qx}a(n)n^{-s} + (\frac{2\pi}{q})^s \frac{1}{\pi}\Gamma(1-s) \cdot \sum_{n \leq y}\sin(\frac{\pi s}{2} + \frac{2\pi na}{q})n^{s-1} + q^{-s} \cdot E(s,x,\frac{a}{q}),$$

where $x > 0$ arbitrary, $y = \frac{|t|}{2\pi x}$; and $a(n) = 1$ for $n \equiv a(\bmod q)$ and 0

otherwise. In view of the fact that $L(s,\chi) = \sum_{a=1}^{q} \chi(a) z(s,a,q)$, this results in

the AFE for $L(s,\chi)$. Following our Theorem 1, this will result in the AFE

for $L(s,\chi_1)L(s,\chi_2)$, where $L(s,\chi_1)$ and $L(s,\chi_2)$ are Dirichlet series

corresponding to Dirichlet characters $\chi_1(\bmod q_1)$ and $\chi_2(\bmod q_2)$

respectively, where $q_1, q_2 \geq 1$ are given integers .. In view of the highly generalised identity stated in the sequel, we can actually give an AFE for the product of functions represented by generalised Dirichlet series, instead of just ordinary Dirichlet series in some half planes.

Next, we state our Theorem 1.

**Theorem 1**: Let the functions $f_1(s)$ and $f_2(s)$ satisfy in a certain region of

complex plane, the relations

$$f_1(s) = \sum_{n \leq x_1} a_1(n) n^{-s} + \psi_1(s) \sum_{n \leq y_1} b_1(n) n^{s-\delta} + E_1(s, x_1) \text{ and}$$

$$f_2(s) = \sum_{n \leq x_2} a_2(n) n^{-s} + \psi_2(s) \sum_{n \leq y_2} b_2(n) n^{s-\delta} + E_2(s, x_2)$$

with $x_1, x_2 > 0$, $y_1 = \frac{T_1}{x_1}$ and $y_2 = \frac{T_2}{x_2}$ so that $y_1, y_2 \geq 0$,

where $\psi_1(s)$ and $\psi_2(s)$ are given functions and $\delta$ is a fixed complex number.

Then in the same region, we have a relation for $f_1(s) \cdot f_2(s)$ namely,

$$f_1(s) f_2(s) = \sum_{n \leq x_1 x_2} A(n) n^{-s} + \psi_1(s) \psi_2(s) \sum_{n \leq y_1 y_2} B(n) n^{s-\delta} + E_{1,2}(s, x_1 x_2),$$

where $A(n) = \sum_{uv=n} a_1(u) a_2(v)$; $B(n) = \sum_{uv=n} b_1(u) b_2(v)$;

and $E_{1,2}(s, x_1 x_2) = \sum_{n \leq x_1} a_1(n) n^{-s} \cdot E_2(s, \frac{x_1 x_2}{n}) + \sum_{n \leq x_2} a_2(n) n^{-s} \cdot E_1(s, \frac{x_1 x_2}{n})$

$+ \psi_1(s) \sum_{m \leq y_1} b_1(m) m^{s-\delta} \cdot E_2(s, \frac{m x_2}{y_1}) + \psi_2(s) \sum_{m \leq y_2} b_2(m) m^{s-\delta} \cdot E_1(s, \frac{m x_1}{y_2})$

$+ E_1(s, x_1) E_2(s, x_2) + L_1 + L_2$

$= I_1 + I_2 + I_3 + I_4 + I_5 + L_1 + L_2$ , say ,

where $L_1 = \psi_1(s) \ell_2^{-s} \cdot \sum_{m \leq y_1} b_1(m) \cdot a_2(\ell_2 m) m^{-\delta}$

and $L_2 = \psi_2(s) \ell_1^{-s} \cdot \sum_{n \leq y_2} b_2(n) \cdot a_1(\ell_1 n) n^{-\delta}$ , with $\ell_1 = \frac{x_1}{y_2}$ and $\ell_2 = \frac{x_2}{y_1}$ .

Here $a_1(\ell_1 n) = 0$ , if $\ell_1 n$ is not an integer ; and

$a_2(\ell_2 m) = 0$ , if $\ell_2 m$ is not an integer .

In addition , if there are functional equations

$f_1(s) = \psi_1(s) f_1(\delta - s)$ and $f_2(s) = \psi_2(s) f_2(\delta - s)$, then we can write

$I_3 = \psi_1(s) \psi_2(s) \sum_{n \leq y_1} b_1(n) n^{s-\delta} E_2(\delta - s, \frac{y_1 y_2}{n})$ and

$I_4 = \psi_1(s) \cdot \psi_2(s) \sum_{n \leq y_2} b_2(n) n^{s-\delta} \cdot E_1(\delta - s, \frac{y_1 y_2}{n})$ .

Here if either $y_1$ or $y_2$ (or both) are zeros or less than 1, then the corresponding summations $\sum_{n \leq y_1}$ or $\sum_{n \leq y_2}$ are empty . We treat empty sums to be zero .

**Remark** : If $\ell_2 = \frac{x_2}{y_1}$ is irrational , then $L_1 = 0$ and if $\ell_1 = \frac{x_1}{y_2}$ is irrational ,

then $L_2 = 0$ .

**Corollary 1** : a) Suppose a function $f(s)$ has an AFE in a certain region as

: 7 :

$$f(s) = \sum_{n\le x} a(n)n^{-s} + \psi(s) \sum_{n\le y} b(n)n^{s-\delta} + E_1(s,x)$$ with $x > 0, y = \frac{T}{x}$.

Then the AFE for $f^2(s)$ is given by

$$f^2(s) = \sum_{n\le x^2} A(n)n^{-s} + \psi^2(s) \sum_{n\le y^2} B(n)n^{s-\delta} + E_2(s,x^2)$$, where

$$E_2(s,x^2) = 2\left\{\sum_{n\le x} a(n)n^{-s} E_1(s,\tfrac{x^2}{n}) + \psi(s)\sum_{n\le y} b(n)n^{s-\delta} E_1(s,\tfrac{nx}{y}) + L\right\} + E_1^2(s,x),$$

where $L = \psi(s)\ell^{-s} \sum_{n\le y} b(n)a(\ell n)n^{-\delta}$ with $a(\ell n) = 0$, if $\ell n$ is not an integer.

Here $\ell = \frac{x}{y}$, $A(n) = \sum_{uv=n} a(u)a(v)$ and $B(n) = \sum_{uv=n} b(u)b(v)$.

In addition, if $f(s)$ has a functional equation $f(s) = \psi(s) f(\delta - s)$, then

$$E_2(s,x^2) = 2\left\{\sum_{n\le x} a(n)n^{-s} \cdot E_1(s,\tfrac{x^2}{n}) + \psi^2(s)\sum_{n\le y} b(n)n^{s-\delta} \cdot E_1(\delta - s,\tfrac{y^2}{n}) + L\right\} + E_1^2(s,x).$$

**Remark**: If $\ell = \frac{x}{y}$ is irrational, then $L = 0$.

Note that if $y = 0$, this corollary 1 reads as follows.

b) If $f(s) = \sum_{n\le x} a(n)n^{-s} + E_1(s,x)$ with $x > 0$,

then $f^2(s) = \sum_{n\le x^2} A(n)n^{-s} + E_2(s,x^2)$,

where $E_2(s,x^2) = 2\sum_{n\le x} a(n)n^{-s} E_1(s,\tfrac{x^2}{n}) + E_1^2(s,x)$ with $A(n) = \sum_{uv=n} a(u)a(v)$.

Thus we get an AFE for $\zeta(s)$ for $s > 0$ and $s \neq 1$ as follows.

For $s > 0$ and $s \neq 1$, let $\zeta(s) = \sum_{n\le x} n^{-s} + E_1(s,x)$,

then $\zeta^2(s) = \sum_{n \leq X} d(n)n^{-s} + E_2(s, X)$,

with $E_2(s, X) = 2 \sum_{n \leq \sqrt{X}} n^{-s} E_1(s, \frac{X}{n}) + E_1^2(s, \sqrt{X})$,

where $x, X > 0$ are arbitrary and d(n) denotes divisor function.

If $\Delta(X)$ is the remainder term of Dirichlet's divisor problem so that

$\sum_{n \leq X} d(n) = X \log X + (2\gamma - 1)X + \frac{1}{4} + \Delta(X)$, $\gamma$ being Euler's constant, it is

well-known that $\Delta(X) = -2 \sum_{n \leq \sqrt{X}} \psi(\frac{X}{n}) + 0(1)$. Here $\psi(u) = u - [u] - \frac{1}{2}$. In the light of the presence of the term $-x^{-s}\psi(x)$ (in $E_1(s, x)$), which contributes the

term $-2X^{-s} \sum_{n \leq X} \psi(\frac{X}{n}) = X^{-s}(\Delta(X) + 0(1))$ to $E_2(s, X)$ (via the term

$2 \sum_{n \leq \sqrt{X}} n^{-s} E_1(s, \frac{X}{n})$ in $E_2(s, X)$) the connection between $\Delta(X)$ and $E_2(s, X)$ is

obvious for $0 < s \neq 1$. On the other hand, when $s$ is complex in the critical

strip, the term $-x^{-s}\psi(x)$ in $E_1(s, x)$ contributes to $E_2(s, X)$

(via $2 \sum_{n \leq \sqrt{X}} n^{-s} E_1(s, \frac{X}{n}) + \chi^2(s) \sum_{n \leq \sqrt{Y}-s} n^{s-1} E_1(1-s, \frac{Y}{n})$), the expression

$-2\{ \sum_{n \leq \sqrt{X}} n^{-s}(\frac{X}{n})\psi(\frac{X}{n}) + \chi^2(s) \sum_{n \leq \sqrt{Y}} n^{s-1}(\frac{Y}{n})^{s-1}\psi(\frac{Y}{n}) \} = X^{-s}\Delta(X) + \chi^2(s)Y^{s-1}\Delta(Y)$,

where $XY = (\frac{t}{2\pi})^2$.

In view of the AFE for $\zeta(s)$ namely,

$\zeta(s) = \sum_{n \leq X} n^{-s} + \chi(s) \sum_{n \leq Y} n^{s-1} + E_1(s, X)$, where $X > 0$, $Y = \frac{|t|}{2\pi X}$ and

: 9 :

$\zeta(s) = \chi(s)\zeta(1-s)$ , with $E_1(s, X) \ll X^{-\sigma} + \frac{X^{1-\sigma}}{\sqrt{t}}$ provided $X \gg 1$ , $Y \gg 1$ ,

we can give a conjecture for the estimates of the remainder term of the AFE for

$\sqrt{\zeta(s)}$ , where the principal branch of logarithm has been considered .

Let $\sqrt{\zeta(s)} = \sum_{n \geq 1} a(n) n^{-s}$ for $\sigma > 1$ . Note that $\sqrt{\zeta(s)} = \sqrt{\chi(s)} \cdot \sqrt{\zeta(1-s)}$ .

Thus , we have an AFE for $\sqrt{\zeta(s)}$ in the form

$$\sqrt{\zeta(s)} = \sum_{n \leq \sqrt{X}} a(n) n^{-s} + \sqrt{\chi(s)} \sum_{n \leq \sqrt{Y}} a(n) n^{s-1} + E_{\frac{1}{2}}(s, \sqrt{X}) ,$$

Thus we state our Results as follows .

**Result 1** : Let $\sqrt{\zeta(s)} = \sum_{n \geq 1} a(n) n^{-s}$ for $\sigma > 1$ .

Let $0 < \sigma < 1$ ; let $x > 0; xy = \sqrt{\frac{t}{2\pi}}$ , where $t \geq 0$ .

Then $\sqrt{\zeta(s)} = \sum_{n \leq x} a(n) n^{-s} + \sqrt{\chi(s)} \sum_{n \leq y} a(n) n^{s-1} + E_{\frac{1}{2}}(s, x)$ ,

where for $X, Y \gg 1$ and $XY = \frac{t}{2\pi}$ , $X^{-\sigma} + \frac{X^{1-\sigma}}{\sqrt{t}} \gg E_1(s, X)$

$$= 2 \left\{ \sum_{n \leq \sqrt{X}} a(n) n^{-s} \cdot E_{\frac{1}{2}}(s, \tfrac{X}{n}) + \chi(s) \sum_{n \leq \sqrt{Y}} a(n) n^{s-1} E_{\frac{1}{2}}(1-s, \tfrac{Y}{n}) \right.$$

$$\left. + \sqrt{\chi(s)} (\sqrt{\tfrac{X}{Y}})^{-s} \sum_{n \leq \sqrt{X}} \frac{a(n) a(n\sqrt{\tfrac{X}{Y}})}{n} \right\} + E_{1/2}^2(s, \sqrt{X}) .$$

Here the third term in the curly bracket (involving summation) vanishes , if

$\sqrt{\frac{X}{Y}}$ is irrational .

More generally , we have the following

**Result 2** : For an integer $q \geq 1$ , let $\chi(\bmod q)$ be a primitive Dirichlet

character and let $L(s, \chi)$ be the corresponding Dirichlet L-series . Let

$$\sqrt{L(s,\chi)} = \sum_{n\geq 1} a(n,\chi) n^{-s} \text{ for } \sigma > 1 \text{ and let } L(s,\chi) = \psi(s,\chi) L(1-s,\overline{\chi}).$$

Let $x > 0$; $xy = \sqrt{\frac{qt}{2\pi}}$, where $t \geq 0$; $0 \leq \sigma \leq 1$. Then

$$\sqrt{L(s,\chi)} = \sum_{n \leq x} a(n,\chi) n^{-s} + \sqrt{\psi(s,\chi)} \sum_{n \leq y} a(n,\overline{\chi}) n^{s-1} + E_{\frac{1}{2}}(s, x, \chi),$$

where for $X, Y \gg 1$ and $XY = \frac{qt}{2\pi}$, $q^{\frac{1}{2}} X^{-\sigma} + \frac{X^{1-\sigma}}{\sqrt{qt}} \gg E_1(s, X, \chi)$

$$= \left( \sum_{n \leq \sqrt{X}} a(n,\chi) n^{-s} \cdot E_{\frac{1}{2}}(s, \tfrac{X}{n}, \chi) + \psi(s,\chi) \sum a(n,\overline{\chi}) n^{s-1} E_{\frac{1}{2}}(1-s, \tfrac{Y}{n}, \overline{\chi}) + L \right)$$

$+ E_{\frac{1}{2}}^2(s, \sqrt{X}, \chi)$, with $L = \sqrt{\psi(s,\chi)} \ell^{-s} \sum_{n \leq \sqrt{Y}} \frac{a(\ell n, \chi) a(n, \overline{\chi})}{n}$ and $\ell = \sqrt{\frac{X}{Y}}$.

**Note** : Here, we have considered principal branches of logarithms of the functions concerned.

Next we state a few more corollaries of Theorem 1, wherein the error terms (of the product of functions) have been estimated by taking only trivial bounds and hence may be improved upon by estimating the error terms more delicately.

Next, we shall state an AFE for $\zeta(s,\alpha)$ with explicit remainder term in a convenient form and consequently, state an AFE for $L(s,\chi)$, where $\chi(\text{mod } q)$ is a primitive Dirichlet character modulo an integer $q \geq 1$. We shall state these as propositions. These results have been stated indirectly in author [3].

**Proposition 1** : Let $s = \sigma + it$ with $t \geq 0$ and $0 < \sigma < 1$ and let $0 < \alpha \leq 1$. Let $0 < \eta < 1$ be a small number. Let $x > 0$ be arbitrary and let $y = \frac{t}{2\pi x}$. Let $\psi(u) = u - [u] - \frac{1}{2}$. Then we have

$$\zeta(s,\alpha) = \sum_{0 \leq n \leq x-\alpha}(n+\alpha)^{-s} + i(2\pi)^{s-1}\Gamma(1-s) \cdot e^{-\frac{\pi i s}{2}} \sum_{1 \leq n \leq y} e^{-2\pi i n \alpha} \cdot n^{s-1} + E_1(s,x,\alpha),$$

where $E_1(s,x,\alpha) - \frac{x^{1-s}}{s-1} = K_1 + K_2$

with $K_1 = x^{-s}\left(\psi(x-\alpha) + \sum_{1 \leq v < y+\eta} \frac{e^{2\pi i v(x-\alpha)}}{2\pi i v}\right) + \sum_{v > y+\eta} \frac{s}{2\pi i v} \int_x^\infty \frac{e^{2\pi i v(u-\alpha)}}{u^{s+1}}du - \sum_{1 \leq v \leq y-\eta} \int_0^x \frac{e^{2\pi i v(u-\alpha)}}{u^s}du$

$-\sum_{v \geq 1} \frac{s}{2\pi i v} \int_x^\infty \frac{e^{-2\pi i v(x-\alpha)}}{u^{s+1}}du$ and $K_2 = \sum_{y<v<y+\eta} \frac{s}{2\pi i v} \int_x^\infty \frac{e^{2\pi i v(u-\alpha)}}{u^{s+1}}du - \sum_{y-\eta<v \leq y} \int_0^x \frac{e^{2\pi i v(u-\alpha)}}{u^s}du$.

In addition, if $y = v_0 \geq 1$ is an integer, we can write

$E_1(s,x,\alpha) - \frac{x^{1-s}}{s-1} = K_1 + K_2$

with $K_1 = x^{-s}\left(\psi(x-\alpha) + \sum_{1 \leq v \leq v_0} \frac{e^{2\pi i v(x-\alpha)}}{2\pi i v}\right) + \sum_{v \geq v_0+1} \frac{s}{2\pi i v} \int_x^\infty \frac{e^{2\pi i v(u-\alpha)}}{u^{s+1}}du$

$-\sum_{1 \leq v \leq v_0-1} \int_0^x \frac{e^{2\pi i v(u-\alpha)}}{u^s}du - \sum_{v \geq 1} \frac{s}{2\pi i v} \int_x^\infty \frac{e^{-2\pi i v(u-\alpha)}}{u^{s+1}}du$ and $K_2 = -\int_0^x \frac{e^{2\pi i v_0(u-\alpha)}}{u^s}du$.

If $x \gg 1$ and $y \gg 1$, then we have $K_1 \ll x^{-\sigma}\log(t+2)$ and $K_2 \ll \frac{x^{1-\sigma}}{\sqrt{t}}$.

**Note**: If $s > 0$ is real so that $y = 0$, we have $\zeta(s,\alpha) = \sum_{0 \leq n \leq x-\alpha}(n+\alpha)^{-s} + E_1(s,x,\alpha)$,

where $E_1(s,x,\alpha) - \frac{x^{1-s}}{s-1} = x^{-s}\psi(x-\alpha) - s\int_x^\infty \frac{\psi(u-\alpha)}{u^{s+1}}du \ll x^{-\sigma}$.

Let $\chi(\mathrm{mod}\, q)$ be a primitive Dirichlet character modulo an integer $q \geq 1$.

Note that $L(s,\chi) = \sum_{a=1}^q \chi(a) \cdot \left(q^{-s}\zeta(s,\frac{a}{q})\right)$. We write $\tau(\chi,n) = \sum_{a=1}^q \chi(a)e^{\frac{2\pi i n a}{q}}$ and note that $\tau(\chi,n) = \overline{\chi}(n)\tau(\chi)$, where $\tau(\chi) = \tau(\chi,1)$ with $|\tau(\chi)| = q^{\frac{1}{2}}$, if $\chi(\mathrm{mod}\, q)$ is primitive. Also note that if $x > 0$ is arbitrary, then $\sum_{a=1}^{q} \chi(a)\left(\frac{x^{1-s}}{s-1}\right) = 0$ for any non-principal character $\chi(\bmod q)$, where $s \neq 1$. For arbitrary $X > 0$, replacing $x$ by $\frac{X}{q}$ and $\alpha$ by $\frac{a}{q}$ in Proposition 1, we get the following Proposition 2.

**Proposition 2**: Let $0 \leq \sigma \leq 1$ and $t \geq 0$. Let $X > 0$ be arbitrary. Let $\chi(\bmod q)$ be a primitive character and let $2\pi \frac{X}{q} y = t$ and let $0 < \eta < 1$ be a small number. Then we have

$$L(s, \chi) = \sum_{n \leq X} \chi(n) n^{-s} + iq^{-s}(2\pi)^{s-1} \Gamma(1-s) \chi(-1) \sum_{n \leq y} \overline{\chi}(n) n^{s-1} + E_1(s, X, \chi),$$

where $E_1(s, X, \chi) = K_1(\chi) + K_2(\chi)$ with

$$K_1(\chi) = X^{-s} \sum_{a=1}^{q} \chi(a) \left( \psi(\tfrac{X-a}{q}) + \sum_{1 \leq v < y+\eta} \frac{e^{2\pi i v (\frac{X-a}{q})}}{2\pi i v} \right) + \tau(\chi) \sum_{v > y+\eta} \frac{s\overline{\chi}(-v)}{2\pi i v} \int_{X/q}^{\infty} \frac{e^{2\pi i v u}}{q^s \cdot u^{s+1}} du$$

$$- \tau(\chi) \sum_{1 \leq v < y-\eta} \overline{\chi}(-v) \int_{0}^{X/q} \frac{e^{2\pi i v u}}{(qu)^s} du - \tau(\chi) \sum_{v \geq 1} \frac{s\overline{\chi}(v)}{2\pi i v} \int_{X/q}^{\infty} \frac{e^{-2\pi i v u}}{q^s \cdot u^{s+1}} du$$

and $K_2(\chi) = \tau(\chi) \sum_{y < v < y+\eta} \frac{s\overline{\chi}(-v)}{2\pi i v} \int_{X/q}^{\infty} \frac{e^{2\pi i v u}}{q^s \cdot u^{s+1}} du - \tau(\chi) \sum_{y-\eta < v \leq y} \overline{\chi}(-v) \int_{0}^{X/q} \frac{e^{2\pi i v u}}{(qu)^s} du$.

In addition if $y = v_0$ is an integer, we have $K_2(\chi) = -\tau(\chi) \overline{\chi}(-v) \int_{0}^{X/q} \frac{e^{2\pi i v_0 u}}{(qu)^s} du$.

If $X \gg 1$ and $y \gg 1$, we have $K_1(\chi) \ll q^{\frac{1}{2}} X^{-\sigma} \log q(t+2)$ and $K_2(\chi) \ll \frac{X^{1-\sigma}}{\sqrt{qt}}$.

**Note**: If $s > 0$ is real so that $y = 0$, we have $L(s, \chi) = \sum_{n \leq X} \chi(n) n^{-s} + E_1(s, X, \chi)$,

where $E_1(s, X, \chi) = X^{-s} \sum_{a=1}^{q} \chi(a) \psi(\tfrac{X-a}{q}) - sq^{-s} \sum_{a=1}^{q} \chi(a) \int_{X/q}^{\infty} \frac{\psi(u - \tfrac{a}{q})}{u^{s+1}} du$

: 13 :

$$\ll q^{\frac{1}{2}}\log(q+2)\cdot X^{-s} + qX^{-s} \ll qX^{-s}.$$

**More corollaries of Theorem 1 :**

**Corollary 2** : For a primitive character $\chi(\bmod q)$, we have an AFE for $L^2(s,\chi)$ as follows namely :

Let $0 \le \sigma \le 1$, $t \ge 0$ and let $X > 0$ and let $Y = (\frac{qt}{2\pi})^2 / X$ and let

$L(s,\chi) = \psi(s,\chi) L(1-s,\overline{\chi})$. Then, we have

$$L^2(s,\chi) = \sum_{n \le X} d(n)\chi(n)n^{-s} + \psi^2(s,\chi) \sum_{n \le Y} d(n)\overline{\chi}(n)n^{s-1} + E_2(s,X,\chi), \text{ where}$$

$$E_2(s,X,\chi) = 2\left( \sum_{n \le \sqrt{X}} \chi(n) n^{-s} E_1(s, \tfrac{X}{n}, \chi) + \psi^2(s,\chi) \sum_{n \le \sqrt{Y}} \overline{\chi}(n) n^{s-1} E_1(1-s, \tfrac{Y}{n}, \overline{\chi}) + L \right)$$

$+ E_1^2(s, \sqrt{X}, \chi)$, where $L = \psi(s,\chi)\ell^{-s} \sum_{n \le \sqrt{Y}} \dfrac{\chi(\ell n)\overline{\chi}(n)}{n}$ with $\ell = \sqrt{\tfrac{X}{Y}}$.

Here $\chi(\ell n) = 0$, if $\ell n$ is not an integer.

In addition, if $X, Y \gg 1$, we have

$$E_2(s,X,\chi) \ll (q^{\frac{1}{2}} X^{\frac{1}{2}} + \tfrac{X}{\sqrt{qt}} + \sqrt{qt} + q) X^{-\sigma} \log q(t+3).$$

This gives $E_2(s, \tfrac{qt}{2\pi}, \chi) \ll q^{\frac{1}{2}} (qt)^{\frac{1}{2}-\sigma} \cdot \log q(t+3)$.

Following our Theorem 1 and the Proposition 2 above, we shall give an explicit expression for $E_2(s, \tfrac{qt}{2\pi} \chi)$ in the form of Theorem 2 below, which we state without proof. (Here it may be noted that $\psi(s,\chi)\psi(1-s,\overline{\chi}) = 1$).

**Theorem 2** : Let $0 \le \sigma \le 1$ and $t \ge 0$. Let $\chi(\bmod q)$ be a primitive character modulo an integer $q \geq 1$ and let $L(s,\chi) = \psi(s,\chi)L(1-s,\overline{\chi})$. Let x, X>0 be arbitrary and let $xy = \frac{qt}{2\pi}$; $XY = (\frac{qt}{2\pi})^2$ and let $\psi(u) = u - [u] - \frac{1}{2}$.

Let $L(s,\chi) = \sum_{n \leq x} \chi(n)n^{-s} + \psi(s,\chi)\sum_{n \leq y} \overline{\chi}(n)n^{s-1} + E_1(s,x,\chi)$ ;

$L^2(s,\chi) = \sum_{n \leq X} d(n)\chi(n)n^{-s} + \psi^2(s,\chi)\sum_{n \leq Y} d(n)\overline{\chi}(n)n^{s-1} + E_2(s,X,\chi)$ and

$M(s,\chi) = \sum_{n \leq \sqrt{\frac{qt}{2\pi}}} \chi(n)n^{-s} E_1(s,\frac{qt}{2\pi n},\chi)$ ,

where $E_1(s,\frac{qt}{2\pi n},\chi) = (\frac{qt}{2\pi n})^{-s}\sum_{a=1}^{q}\chi(a)\left(\psi(\frac{t}{2\pi n} - \frac{a}{q}) + \sum_{|v| \leq n}\frac{e^{2\pi i v(\frac{t}{2\pi n} - \frac{a}{q})}}{2\pi i v}\right)$

$+ q^{-s}\tau(\chi)\sum_{|v| \geq n+1}\frac{\overline{\chi}(-v)}{2\pi i v}\int_{t/2\pi n}^{\infty}\frac{se^{2\pi i vu}}{u^{s+1}}du - q^{-s}\tau(\chi)\sum_{1 \leq |v| \leq n}\overline{\chi}(-v)\int_{0}^{t/2\pi n}\frac{e^{2\pi i vu}}{u^s}du$ .

Then we have $\psi(1-s,\overline{\chi})E_2(s,\frac{qt}{2\pi},\chi)$

$= 2\left\{\psi(1-s,\overline{\chi})M(s,\chi) + \psi(s,\chi)M(1-s,\overline{\chi}) + \sum_{n \leq \sqrt{\frac{qt}{2\pi}}}'\frac{1}{n}\right\} + \psi(1-s,\overline{\chi})E_1^2(s,\sqrt{\frac{qt}{2\pi}},\chi)$ ,

with the dash over $\sum_{n \leq \sqrt{\frac{qt}{2\pi}}}$ denoting the restriction $(n,q) = 1$ , and

$\psi(1-s,\overline{\chi})E_1^2(s,\sqrt{\frac{qt}{2\pi}},\chi) << q^{\frac{1}{2}}t^{-\frac{1}{2}}$ .

**Remark** : It is to be noted that Motohashi [2] has proved (in the case $q = 1$, that is, in the case of Riemann zeta function) that

$\chi(1-s)E_2(s,\frac{t}{2\pi}) = -\sqrt{2}(\frac{t}{2\pi})^{-\frac{1}{2}}\Delta(\frac{t}{2\pi}) + 0(t^{-\frac{1}{4}})$ ,

following Riemann-Siegel formula for $\zeta(s)$. On the same lines, one may obtain results for $\psi(1-s,\overline{\chi})E_2(s,\frac{qt}{2\pi},\chi)$ using Riemann-Siegel type formula for $\zeta(s,\alpha)$ of author [5]. Of course, the Riemann-Siegel type formula for

$\zeta(s,\alpha)$ could be useful, while dealing with $E_2(s,\frac{qt}{2\pi},\chi)$ as a function of t.

However, while dealing with $E_2(s,\frac{qt}{2\pi},\chi)$ as a function of q, the expression for

$E_2(s,\frac{qt}{2\pi},\chi)$ of our Theorem 2 above, is useful.

**Additional Corollaries of Theorem 1**:

**Corollary 3** : I) If $X > (\frac{c}{2\pi}t)^2$ for an absolute constant $c > 1$,

we have for $s \neq 1$,

$$\zeta^2(s) = \sum_{n \leq X} \frac{d(n)}{n^s} + \frac{2X^{1-s}}{(s-1)} \cdot \sum_{n \leq \sqrt{X}} \frac{1}{n} + 0(X^{\frac{1}{2}-\sigma}).$$

II) If for integer $q \geq 1$, $\chi(\mod q)$ is a primitive character, then for $X > (\frac{c}{2\pi}qt)^2$

with constant $c > 1$ and for $t \geq 1$, we have

$$L^2(s,\chi) = \sum_{n \leq X} \frac{d(n)\chi(n)}{n^s} + 0(q^{1/2}X^{1/2-\sigma})$$

**Remark** : I) of Corollary 3 follows from $\zeta(s) = \sum_{n \leq x} n^{-s} + \frac{x^{1-s}}{s-1} + 0(x^{-\sigma})$

for $x > \frac{ct}{2\pi}$ and for $0 < \sigma < 1$.

II) of Corollary 3 follows from $L(s,\chi) = \sum_{n \leq x} \frac{\chi(n)}{n^s} + 0(q^{1/2}x^{-\sigma}\log(q+2))$

for $x > \frac{c}{2\pi}qt$ and $0 \leq \sigma \leq 1$. ( See corollary of author [3] ).

**Corollary 4** : For integers $q_1, q_2 \geq 1$, let $\chi_1(\mod q_1)$ and $\chi_2(\mod q_2)$ be Dirichlet

characters both primitive and let $L(s,\chi_1)$ and $L(s,\chi_2)$ be the corresponding

Dirichlet L-series. Let $x_1, y_1, x_2, y_2 \gg 1$ with $2\pi x_1 y_1 = q_1 t$ and $2\pi x_2 y_2 = q_2 t$ and

Let $L(s,\chi_i) = \psi(s,\chi_i)L(1-s,\overline{\chi_i})$ for $i = 1,2$.

Then we have

$$L(s,\chi_1)L(s,\chi_2) = \sum\sum_{mn \le x_1 x_2}\chi_1(m)\chi_2(n)(mn)^{-s} + \psi(s,\chi_1)\psi(s,\chi_2)\sum\sum_{mn \le y_1 \cdot y_2}\overline{\chi_1}(m)\overline{\chi_2}(n)(mn)^{s-1}$$
$$+ E(s,x_1 x_2,\chi_1,\chi_2), \text{ where } E(s,x_1 x_2,\chi_1,\chi_2)$$

$$\ll (x_1 x_2)^{-\sigma} \log q_1 q_2 (t+2) \cdot$$
$$\cdot \{x_1\sqrt{q_2} + x_2\sqrt{q_1} + x_1 x_2 t^{-\frac{1}{2}}(q_1^{-\frac{1}{2}} + q_2^{-\frac{1}{2}}) + \sqrt{q_1 t} + \sqrt{q_2 t} + \sqrt{q_1 q_2}\}.$$

Taking $x_1 = \sqrt{\frac{q_1 t}{2\pi}}$ and $x_2 = \sqrt{\frac{q_2 t}{2\pi}}$, we have

$$E(s,\sqrt{q_1 q_2}\tfrac{t}{2\pi},\chi_1,\chi_2) \ll (q_1 q_2)^{\frac{1}{2}-\frac{\sigma}{2}} t^{\frac{1}{2}-\sigma} \log q_1 q_2 (t+2).$$

From Corollary 2, we get the following corollary.

**Corollary 5**: If $\chi(\bmod q)$ is a primitive Dirichlet character and $L(s,\chi)$ is the corresponding Dirichlet L-series, then for

$x_1, x_2 \gg 1; y_1, y_2 \gg 1; 2\pi x_1 y_1 = q\,t; 2\pi x_2 y_2 = t$, we have

$$\zeta(s)L(s,\chi) = \sum_{n \le x_1 x_2} A(n)n^{-s} + \chi(s)\psi(s,\chi)\sum_{n \le y_1 y_2}\overline{A(n)}n^{s-1} + E(s,x_1 x_2,\chi),$$

where $A(n) = \sum_{d|n}\chi(d); \zeta(s) = \chi(s)\zeta(1-s); L(s,\chi) = \psi(s,\chi)L(1-s),\overline{\chi})$

and $E(s,x_1 x_2,\chi) \ll (\sqrt{q}\cdot x_1 + x_2 + \frac{x_1 x_2}{\sqrt{t}} + \sqrt{qt})(x_1 x_2)^{-\sigma} \cdot \log q(t+2)$.

For integral $r \ge 1$ and for $0 < \sigma < 1$, consider the AFE of $\zeta^r(s)$ in the

form, $\zeta^r(s) = \sum_{n \le X} d_r(n)n^{-s} + \chi^r(s)\sum_{n \le Y} d_r(n)n^{s-1} + E_r(s,X)$, where $XY = (\frac{t}{2\pi})^r$;

For $X \gg 1$ ; $Y \gg 1$, assuming the well-known bounds (from literature),

$E_1(s,X) \ll (X^{-\sigma} + \frac{X^{1-\sigma}}{\sqrt{t}})\log(|t|+2)$ and $E_2(s,X) \ll X^{\frac{1}{2}-\sigma}\log(|t|+2)$, we get the

: 17 :

following bounds for $E_3(s,X)$ and $E_4(s,X)$, which we state as Corollary 6.

**Corollary 6** : I) Let $0 \leq \sigma < 1$; $t \geq 0$; and $XY = (\frac{t}{2\pi})^3$, where $X >> 1$; $Y >> 1$.

Then $E_3(s, X) << (X^{\frac{2}{3}} + \frac{X}{\sqrt{t}} + \sqrt{t} + t^{1-\sigma})X^{-\sigma} \log^2(t+2)$.

II) Let $0 \leq \sigma < 1$; $t \geq 0$; and $\epsilon > 0$ be arbitrary. Let $XY = (\frac{t}{2\pi})^4$, where $X >> 1$; $Y >> 1$. Then $E_4(s, X) << (X^{\frac{3}{4}} + tX^{\frac{1}{4}})X^{-\sigma} \cdot t^{\epsilon}$.

From Corollary 2, putting $q = 1$, we get

$E_2(s, X) << (\frac{X}{\sqrt{t}} + \sqrt{t})X^{-\sigma} \log(t+3)$, so that $E_2(s, \frac{t}{2\pi}) << t^{\frac{1}{2}-\sigma} \log(t+3)$.

In fact, a slightly superior inequality namely, $E_2(s, X) << X^{\frac{1}{2}-\sigma} \log(t+3)$, is already known giving $E_2(s, \frac{t}{2\pi}) << t^{\frac{1}{2}-\sigma} \log(t+3)$, which coincides with our estimate. Though our estimate is somewhat inferior, we have to appreciate the fact that our estimate for $E_2(s, X)$ makes use of bare minimal information namely $E_1(s, x) << x^{-\sigma} \log(y+3) + \frac{x^{1-\sigma}}{\sqrt{t}}$ and nothing else; while the proof of $E_2(s, X) << X^{\frac{1}{2}-\sigma} \log(t+3)$ makes use of a good deal of information. Incidentally, Motohashi[1] has proved better namely $E_2(s, \frac{t}{2\pi}) << t^{\frac{1}{3}-\sigma} \log(t+3)$.

It is not difficult to see that we can also give an AFE for $f_1(s)f_2(\bar{s})$ in the form, $f_1(s)f_2(\bar{s}) = \sum\sum_{mn \leq x_1 x_2} \frac{a_1(m)m^{-it}a_1(n)n^{it}}{(mn)^\sigma}$

$+ \psi_1(s)\psi_2(\bar{s})\sum\sum_{mn \leq y_1 y_2} b_1(m)m^{it}b_2(n)n^{-it}(mn)^{\sigma-\delta} + H_{1,2}(s, x_1 x_2)$

for $\delta$ real, with an explicit expression for $H(s, x_1 x_2)$ in terms of $E_1(s, x_1)$ and

: 18 :

$E_2(\bar{s}, x_2)$. Also note that the remainder term $H_{1,2}(s, x_1 x_2)$ can be estimated, if we

have estimates for $E(s, x_1)$, as $x_1$ varies; and $E(s, x_2)$, as $x_2$ varies; and have estimates for arithmetical functions $a_1(n), a_2(n), b_1(n)$ and $b_2(n)$. Applying the above Theorem 1, r times (for integral $r \geq 1$) to the function $f_1(s)$, we get

$$f_1^r(s) = \sum_{n_1 n_2 \ldots n_r \leq x_1^r} \frac{a_1(n_1) a_1(n_2) \ldots a_1(n_r)}{(n_1 n_2 \ldots n_r)^s} + \psi_1^r(s) \sum_{n_1 n_2 \ldots n_r \leq y_1^r} b_1(n_1) b_1(n_2) \ldots b_1(n_r)(n_1 n_2 \ldots n_r)^{s-\delta}$$

+ remainder term. Thus for any integer $r \geq 1$, if $\zeta^{\frac{1}{r}}(s) = \sum_{n \geq 1} a_r(n) n^{-s}$ for $\sigma > 1$;

where $\zeta(s)$ is Riemann zeta function, in view of the functional equation

$\zeta^{\frac{1}{r}}(s) = \chi^{\frac{1}{r}}(s) \zeta^{\frac{1}{r}}(1-s)$ with principal branches of logarithm under consideration,

we have an AFE for $\zeta^{\frac{1}{r}}(s)$ namely,

$$\zeta^{\frac{1}{r}}(s) = \sum_{n \leq X} a_r(n) n^{-s} + \chi^{\frac{1}{r}}(s) \sum_{n \leq Y} a_r(n) n^{s-1} + E_{\frac{1}{r}}(s, X) \text{ with } XY = (\tfrac{t}{2\pi})^{\frac{1}{r}}.$$

It may be possible to estimate $E_{\frac{1}{r}}(s, X)$ for $t >> 1$, in view of the estimate

$E_1(s, x) << x^{-\sigma} \log(y+2) + x^{1-\sigma} t^{-\frac{1}{2}}$ for $t >> 1$ where $E_1(s, x)$ is the remainder term of AFE for $\zeta(s)$.

The main idea in the proof of our Theorem is a very simple identity namely, for any function $\phi(x, y)$ of two variables and for real $x_1, x_2 \geq 0$, we have

I) $\sum \sum_{mn \leq x_1 x_2} \phi(m,n) = \sum_{m \leq x_1} \sum_{n \leq \frac{x_1 x_2}{m}} \phi(m,n) + \sum_{n \leq x_2} \sum_{m \leq \frac{x_1 x_2}{n}} \phi(m,n) - \sum_{m \leq x_1} \sum_{n \leq x_2} \phi(m,n)$

or equivalently $\sum_{m \leq x_1} \sum_{n \leq x_2} \phi(m,n) = \sum \sum_{mn \leq x_1 x_2} \phi(m,n) - \sum_{m \leq x_1} \sum_{x_2 < n \leq \frac{x_1 x_2}{m}} \phi(m,n)$

$- \sum_{n \leq x_2} \sum_{x_1 < m \leq \frac{x_1 x_2}{n}} \phi(m,n)$.

: 19 :

More generally, instead of considering sets of lattice points $(m,n)$, we can generalise this identity to set of points of the form $(\lambda_m, \lambda_n)$, where $(\lambda_n)_{n=1}^{\infty}$ is a strictly increasing sequence of positive real numbers such that $\lambda_n \to \infty$ as $n \to \infty$. The more general type of identity will have to be used, while considering AFEs for functions defined via generalised Dirichlet Series of the form $\sum_n a_n \lambda_n^{-s}$, instead of functions defined via Dirichlet Series of the form $\sum_{n \geq 1} a_n n^{-s}$. The more general identity will read as follows, namely,

II) $\sum_{\lambda_m \leq x_1} \sum_{\lambda_n \leq x_2} \phi(\lambda_m, \lambda_n) = \sum_{\lambda_m, \lambda_n \leq x_1 x_2} \sum \phi(\lambda_m, \lambda_n)$

$- \sum_{\lambda_m \leq x_1} \sum_{x_2 < \lambda_n \leq \frac{x_1 x_2}{\lambda_m}} \phi(\lambda_m, \lambda_n) - \sum_{\lambda_n \leq x_2} \sum_{x_1 < \lambda_m \leq \frac{x_1 x_2}{\lambda_n}} \phi(\lambda_m, \lambda_n)$

Here it must be mentioned that following Dirichlet's device, Motohashi [1] obtained an identity namely, for any integer $N \geq 1$,

$\sum_{n \leq N} d(n) a(n) = 2 \sum_{n \leq \sqrt{N}} \sum_{m \leq \frac{N}{n}} a(mn) - \sum_{m, n \leq \sqrt{N}} \sum a(mn)$

for any arithmetical function $a(n)$. Here $d(n)$ is divisor function.

Our identity above can be easily verified by considering points $(\lambda_m, \lambda_n)$ between $X - Y$ axes and the rectangular hyperbola $XY = x_1 x_2$ and noting that this region can be divided into three obvious disjoint regions. Motohashi's identity is a particular case of our identity for $\lambda_n = n$, $\phi(m, n) = a(mn)$ and $x_1 = x_2 = \sqrt{N}$; and Dirichlet's device, a particular case of Motohashi's identity for constant function $a(n) = 1$.

: 20 :

Our motivation comes from Motohashi [1], who uses his identity to obtain an expression $E_2(s, X)$, the remainder term of AFE for $\zeta^2(s)$ in the particular case $X = \frac{t}{2\pi}$. Motohashi [1] derives an expression $E_2(s, X)$

$$= 2\left\{ \sum_{n \leq \sqrt{X}} n^{-s} E(s, \tfrac{X}{n}) + \chi^2(s) \sum_{n \leq \sqrt{X}} n^{s-1} E(1-s, \tfrac{X}{n}) + \chi(s) \sum_{n \leq \sqrt{X}} \tfrac{1}{n} \right\} + E^2(s, \sqrt{X})$$

in the particular case of $X = \frac{t}{2\pi}$, where $E(s, x)$ is the remainder term of AFE for $\zeta(s)$. Following Riemann-Siegel's method, Motohashi [1] uses his expression for $E_2(s, \frac{t}{2\pi})$ to obtain a counterpart of Riemann-Siegel's formula for $\zeta^2(s)$ and also the estimate $E_2(s, \frac{t}{2\pi}) \ll t^{\frac{1}{3} - \sigma} \log t$ for $t \geq 3$. In particular, he obtains the exact relation between $E_2(s, \frac{t}{2\pi})$ and $\Delta(\frac{t}{2\pi})$, where $\Delta(X)$ is the remainder term of Dirichlet's divisor problem. Motohashi [2] also obtains an expression for $E_2(s, \frac{\alpha t}{2\pi})$ for a rational $\alpha$. However, his expressions for $E_2(s, \frac{\alpha t}{2\pi})$ are not as precise as in the case of $E_2(s, \frac{t}{2\pi})$.

**Proof of Theorem 1 :** We have in the region under consideration for $x_1, x_2 > 0$ and $y_1, y_2 \geq 0$,

$$f_1(s) = \sum_{m \leq x_1} a_1(m) m^{-s} + \psi_1(s) \sum_{m \leq y_1} b_1(m) m^{s-\delta} + E_1(s, x_1) \text{ and}$$

$$f_2(s) = \sum_{n \leq x_2} a_2(n) n^{-s} + \psi_2(s) \sum_{n \leq y_2} b_2(n) n^{s-\delta} + E_2(s, x_2);$$

Note that if $y_1 < 1$ (or $y_2 < 1$), then the corresponding summation $\sum_{n \leq y_1}$ is empty

: 21 :

(or $\sum_{n \leq y_2}$ is empty). Also note that $A(r) = \sum_{mn=r} a_1(m)a_2(n)$.

Thus $\sum_{r \leq x_1 x_2} \frac{A(r)}{r^s} = \sum_{r \leq x_1 x_2} \frac{(\sum_{mn=r} a_1(m)a_2(n))}{r^s} = \sum_{r \leq x_1 x_2} \frac{(\sum_{mn=r} a_1(m)a_2(n))}{(mn)^s} = \sum_{r \leq x_1 x_2} \sum_{mn=r} \frac{a_1(m)}{m^s} \cdot \frac{a_2(n)}{n^s}$

$= \sum_{mn \leq x_1 x_2} \sum \frac{a_1(m)}{m^s} \cdot \frac{a_2(n)}{n^s}$. Similarly $\sum_{r \leq y_1 y_2} B(r) r^{s-\delta} = \sum_{mn \leq y_1 y_2} \sum b_1(m) m^{s-\delta} \cdot b_2(n) n^{s-\delta}$.

Thus

$$f_1(s)f_2(s) = \sum_{m \leq x_1} a_1(m)m^{-s} \sum_{n \leq x_2} a_2(n)n^{-s}$$

$$+ \psi_1(s)\psi_2(s) \sum_{m \leq y_1} b_1(m)m^{s-\delta} \sum_{n \leq y_2} b_2(n)n^{s-\delta} + E_1(s,x_1)E_2(s,x_2)$$

$$+ \psi_2(s) \sum_{m \leq x_1} a_1(m)m^{-s} \sum_{n \leq y_2} b_2(n)n^{s-\delta} + \psi_1(s) \sum_{m \leq x_2} a_2(m)m^{-s} \sum_{n \leq y_1} b_1(n)n^{s-\delta}$$

$$+ E_1(s,x_1) \left( \sum_{n \leq x_2} a_2(n)n^{-s} + \psi_2(s) \sum_{n \leq y_2} b_2(n)n^{s-\delta} \right) + E_2(s,x_2) \left( \sum_{n \leq x_1} a_1(n)n^{-s} + \psi_1(s) \sum_{n \leq y_1} b_1(n)n^{s-\delta} \right)$$

Next

$$\sum_{m \leq x_1} a_1(m)m^{-s} \sum_{n \leq x_2} a_2(n)n^{-s} = \sum_{n \leq x_1 x_2} A(n)n^{-s}$$

$$- \left( \sum_{m \leq x_1} a_1(m)m^{-s} \sum_{x_2 < n \leq \frac{x_1 x_2}{m}} a_2(n)n^{-s} + \sum_{n \leq x_2} a_2(n)n^{-s} \sum_{x_1 < m \leq \frac{x_1 x_2}{n}} a_1(m)m^{-s} \right)$$

Similarly

$$\psi_1(s)\psi_2(s) \sum_{m \leq y_1} b_1(m)m^{s-\delta} \sum_{n \leq y_2} b_2(n)n^{s-\delta} = \psi_1(s)\psi_2(s) \sum_{n \leq y_1 y_2} B(n)n^{s-\delta}$$

$$- \psi_1(s)\psi_2(s) \left( \sum_{m \leq y_1} b_1(m)m^{s-\delta} \sum_{y_2 < n \leq \frac{y_1 y_2}{m}} b_2(n)n^{s-\delta} + \sum_{n \leq y_2} b_2(n)n^{s-\delta} \sum_{y_1 < m \leq \frac{y_1 y_2}{n}} b_1(m)m^{s-\delta} \right)$$

: 22 :

Thus, $f_1(s)f_2(s) = \sum_{n \leq x_1 x_2} A(n)n^{-s} + \psi_1(s)\psi_2(s) \sum_{n \leq y_1 y_2} B(n)n^{s-\delta}$

$- \Big\{ \sum_{m \leq x_1} a_1(m)m^{-s} \sum_{x_2 < n \leq \frac{x_1 x_2}{m}} a_2(n)n^{-s} + \sum_{n \leq x_2} a_2(n)n^{-s} \sum_{x_1 < m \leq \frac{x_1 x_2}{n}} a_1(m)m^{-s}$

$+ \psi_1(s)\psi_2(s) \sum_{m \leq y_1} b_1(m)m^{s-\delta} \sum_{y_2 < n \leq \frac{y_1 y_2}{m}} b_2(n)n^{s-\delta}$

$+ \psi_1(s)\psi_2(s) \sum_{n \leq y_2} b_2(n)n^{s-\delta} \sum_{y_1 < m \leq \frac{y_1 y_2}{n}} b_1(m)m^{s-\delta}$

$- \psi_2(s) \sum_{m \leq x_1} a_1(m)m^{-s} \sum_{n \leq y_2} b_2(n)n^{s-\delta}$

$- \psi_1(s) \sum_{m \leq x_2} a_2(m)m^{s-\delta} \sum_{n \leq y_1} b_1(n)n^{s-\delta} - E_1(s, x_1) \Big( \sum_{m \leq x_2} a_2(n)n^{-s} + \psi_2(s) \sum_{n \leq y_2} b_2(n)n^{s-\delta} \Big)$

$- E_2(s, x_2) \Big( \sum_{n \leq x_1} a_1(n)n^{-s} + \psi_1(s) \sum_{n \leq y_1} b_1(n)n^{s-\delta} \Big) - E_1(s, x_1) E_2(s, x_2) \Big\}$

$= \sum_{n \leq x_1 x_2} A(n)n^{-s} + \psi_1(s)\psi_2(s) \sum_{n \leq y_1 y_2} B(n)n^{s-\delta} + E_{1,2}(s, x_1 x_2)$.

Substituting

I) $\sum_{x_2 < n \leq \frac{x_1 x_2}{m}} a_2(n)n^{-s} = \psi_2(s) \sum_{\frac{my_2}{x_1} < n \leq y_2} b_2(n)n^{s-\delta} + E_2(s, x_2) - E_2(s, \frac{x_1 x_2}{m})$

II) $\sum_{x_1 < m \leq \frac{x_1 x_2}{n}} a_1(m)m^{-s} = \psi_1(s) \sum_{\frac{ny_1}{x_2} < m \leq y_1} b_1(m)m^{s-\delta} + E_1(s, x_1) - E_1(s, \frac{x_1 x_2}{n})$

III) $\psi_2(s) \sum_{y_2 < n \leq \frac{y_1 y_2}{m}} b_2(n)n^{s-\delta} = \sum_{\frac{mx_2}{y_1} < n \leq x_2} a_2(n)n^{-s} + E_2(s, x_2) - E_2(s, \frac{mx_2}{y_1})$

IV) $\psi_1(s) \sum_{y_1 < m \leq \frac{y_1 y_2}{n}} b_1(m)m^{s-\delta} = \sum_{\frac{nx_1}{y_2} < m \leq x_1} a_1(m)m^{-s} + E_1(s, x_1) - E_1(s, \frac{nx_1}{y_2})$,

we have the expression in the curly bracket namely, $-E_{1,2}(s, x_1 x_2)$

: 23 :

$$= \psi_2(s) \sum_{m \leq x_1} \frac{a_1(m)}{m^s} \sum_{n \leq y_2} b_2(n) n^{s-\delta} - \psi_2(s) \sum_{m \leq x_1} \frac{a_1(m)}{m^s} \sum_{n \leq \frac{my_2}{x_1}} b_2(n) n^{s-\delta}$$

$$+ E_2(s, x_2) \sum_{m \leq x_1} \frac{a_1(m)}{m^s} - \sum_{m \leq x_1} \frac{a_1(m)}{m^s} E_2(s, \tfrac{x_1 x_2}{m})$$

$$+ \psi_1(s) \sum_{n \leq x_2} \frac{a_2(n)}{n^s} \sum_{m \leq y_1} b_1(m) m^{s-\delta} - \psi_1(s) \sum_{n \leq x_2} \frac{a_2(n)}{n^s} \sum_{m \leq \frac{ny_1}{x_2}} b_1(m) m^{s-\delta}$$

$$+ E_1(s, x_1) \sum_{n \leq x_2} \frac{a_2(n)}{n^s} - \sum_{n \leq x_2} \frac{a_2(n)}{n^s} E_1(s, \tfrac{x_1 x_2}{n}) + \psi_1(s) \sum_{m \leq y_1} b_1(m) m^{s-\delta} \sum_{n \leq x_2} \frac{a_2(n)}{n^s}$$

$$- \psi_1(s) \sum_{m \leq y_1} b_1(m) m^{s-\delta} \sum_{n \leq \frac{mx_2}{y_1}} \frac{a_2(n)}{n^s} - \psi_1(s) \sum_{m \leq y_1} b_1(m) m^{s-\delta} E_2(s, \tfrac{mx_2}{y_1})$$

$$+ E_2(s, x_2) \, \psi_1(s) \sum_{n \leq y_1} b_1(n) n^{s-\delta}$$

$$+ \psi_2(s) \sum_{n \leq y_2} b_2(n) n^{s-\delta} \sum_{m \leq x_1} \frac{a_1(m)}{m^s} - \psi_2(s) \sum_{n \leq y_2} b_2(n) n^{s-\delta} \sum_{m \leq \frac{nx_1}{y_2}} \frac{a_1(m)}{m^s}$$

$$+ \psi_2(s) \, E_1(s, x_1) \sum_{n \leq y_2} b_2(n) n^{s-\delta} - \psi_2(s) \sum_{n \leq y_2} b_2(n) n^{s-\delta} E_1(s, \tfrac{nx_1}{y_2})$$

$$- \psi_2(s) \sum_{m \leq x_1} \frac{a_1(m)}{m^s} \sum_{n \leq y_2} b_2(n) n^{s-\delta} - \psi_1(s) \sum_{m \leq x_2} \frac{a_2(m)}{m^s} \sum_{n \leq y_1} b_1(n) n^{s-\delta}$$

$$= \psi_2(s) \, \Big( \sum_{n \leq y_2} \sum_{m \leq x_1} - \sum_{m \leq x_1} \sum_{n \leq \frac{my_2}{x_1}} - \sum_{n \leq y_2} \sum_{m \leq \frac{nx_1}{y_2}} \Big) b_2(n) n^{s-\delta} a_1(m) m^{-s}$$

$$+ \psi_1(s) \, \Big( \sum_{m \leq y_1} \sum_{n \leq x_2} - \sum_{n \leq x_2} \sum_{m \leq \frac{ny_1}{x_2}} - \sum_{n \leq y_1} \sum_{n \leq \frac{mx_2}{y_1}} \Big) b_1(m) m^{s-\delta} a_2(n) n^{-s}$$

$$- E_1(s, x_1) E_2(s, x_2) - \sum_{m \leq x_1} a_1(m) m^{-s} E_2(s, \tfrac{x_1 x_2}{m})$$

$$-\sum_{n\le x_2}a_2(n)n^{-s}E_1(s,\tfrac{x_1x_2}{n})-\psi_1(s)\sum_{m\le y_1}b_1(m)m^{s-\delta}E_2(s,\tfrac{mx_2}{y_1})-\psi_2(s)\sum_{n\le y_2}b_2(n)n^{s-\delta}E_1(s,\tfrac{nx_1}{y_2})$$

$$=-\psi_2(s)\sum_{\substack{n\le y_2\\m=nx_1/y_2}}b_2(n)n^{s-\delta}a_1(m)m^{-s}-\psi_1(s)\sum_{\substack{m\le y_1\\n=mx_2/y_1}}b_1(m)m^{s-\delta}a_2(n)n^{-s}$$

$$-E_1(s,x_1)E_2(s,x_2)-\sum_{m\le x_1}a_1(m)m^{-s}E_2(s,\tfrac{x_1x_2}{m})$$

$$-\sum_{n\le x_2}a_2(n)n^{-s}E_1(s,\tfrac{x_1x_2}{n})-\psi_1(s)\sum_{m\le y_1}b_1(m)m^{s-\delta}E_2(s,\tfrac{mx_2}{y_1})-\psi_2(s)\sum_{n\le y_2}b_2(n)n^{s-\delta}E_1(s,\tfrac{nx_1}{y_2})$$

This gives the expression for $E_{1,2}(s,x_1x_2)$ as in the statement of our Theorem 1.

**References** :

[1] Y. Motohashi, A note on the approximate functional equation for $\zeta^2(s)$,

   Proc. Japan Acad. Ser 59 A (1983), 392-396.

[2] Y. Motohashi, A note on the approximate functional equation for $\zeta^2(s)$-III,

   ibid 62 A (1986), 410-412

[3] V.V. Rane, On an Approximate Functional Equation for Dirichlet L-Series,

   Math. Ann. 264, 137-145 (1983)

[4] V.V. Rane, A new approximate functional equation for Hurwitz zeta

   function for rational parameter, Proc. Indian acad. Sci. (Math.Sci.), Vol.

   107, No. 4, November 1997, pp 377-385.

[5] V.V. Rane, On the mean square value of Dirichlet L-Series, J. London

   Math. Soc. (2), 21 (1980), 203-215.